# LINEARIZATION COEFFICIENTS FOR ORTHOGONAL POLYNOMIALS USING STOCHASTIC PROCESSES[1]

By Michael Anshelevich

*University of California, Riverside*

Given a basis for a polynomial ring, the coefficients in the expansion of a product of some of its elements in terms of this basis are called linearization coefficients. These coefficients have combinatorial significance for many classical families of orthogonal polynomials. Starting with a stochastic process and using the stochastic measures machinery introduced by Rota and Wallstrom, we calculate and give an interpretation of linearization coefficients for a number of polynomial families. The processes involved may have independent, freely independent or $q$-independent increments. The use of noncommutative stochastic processes extends the range of applications significantly, allowing us to treat Hermite, Charlier, Chebyshev, free Charlier and Rogers and continuous big $q$-Hermite polynomials.

We also show that the $q$-Poisson process is a Markov process.

**1. Introduction.** Let $\{P_n\}$ be a family of polynomials orthogonal with respect to a measure $\mu$ on the real line. One standard combinatorial question is to calculate the moments of the measure, $m_n = \langle x^n \rangle$, where we denote by $\langle \cdot \rangle$ the integral (expectation) with respect to $\mu$. For many classical families of polynomials these moments are positive integers or, more generally, polynomials in parameters with positive integer coefficients. These coefficients beg a combinatorial interpretation, and there exists a large body of work to this effect.

A more general question one can ask is to calculate the linearization coefficients. That is, for $(n_1, n_2, \ldots, n_k)$, we are interested in the expectations $\langle P_{n_1} P_{n_2} \ldots P_{n_k} \rangle$. The name stems from the fact that these are the coefficients

Received February 2003; revised December 2003.
[1]Supported in part by an NSF postdoctoral fellowship.

in the expansion of products of this type in the basis $\{P_n\}$, that is, expansions as sums of orthogonal polynomials. Again, many of these coefficients are positive integers, and so they "count something."

A combinatorial approach to this problem is to construct explicit bijections between structures counted by the linearization coefficients and structures of known cardinality; see, for example, [10]. In this paper we take a different route, and consider a probabilistic interpretation of certain coefficients. The connection to combinatorics is provided by the fact that the moments of a measure are sums, over all set partitions, of products of cumulants of that measure. We will see that certain linearization coefficients can by described in a similar way. The machinery we use is that of stochastic measures, first introduced by Rota and Wallstrom in [12]. In a number of previous papers we extended this machinery from the usual to the noncommutative stochastic processes. This extends the number of polynomial families that we can handle, and so we not only obtain a nice interpretation of known results, but some new results as well. In particular, we show that the linearization coefficients for the continuous big $q$-Hermite polynomials ([11], 3.18) are based on the number of the inhomogeneous set partitions, with an extra statistic counting the number of "restricted crossings" of such partitions.

To be more specific, for each family of polynomials in this paper and the related family of measures of orthogonality, we introduce a, possibly noncommutative, stochastic process $\{X(t)\}$. Then for this process, we introduce a further family $\{\psi_k(t)\}$ of other stochastic processes, which we call full stochastic measures. These objects are orthogonal, and have clean linearization formulas. On the rare occasions when these objects are polynomials in the original process $X$, these formulas translate into the linearization formulas for polynomials.

Another property, which always holds for the full stochastic measures and which in these cases is shared by the orthogonal polynomials, is the martingale property. The Markov property for the $q$-Brownian motion was shown in [6] using the Gaussian properties of the process. Using the above fact for the stochastic measures, we show that the $q$-Poisson process is also a Markov process.

The paper is organized as follows. In Section 2 we describe general combinatorial properties of combinatorial stochastic measures. Section 3 is based on the results of [12] about processes with independent increments, and gives the linearization coefficients for the Hermite and Charlier polynomials. Section 4 is based on the results of [1, 3] about processes with freely independent increments, and gives the linearization coefficients for the Chebyshev polynomials of the 2nd kind and the free Charlier polynomials. Section 5 is based on the results of [2] about $q$-Lévy processes, and gives the linearization coefficients for the continuous and continuous big $q$-Hermite polynomials. It also requires some new results about the $q$-Poisson process. The proofs of these results are contained in the Appendix, which also contains the proof of the Markov property for the $q$-Poisson process.

**2. Combinatorial stochastic measures.** Let $(\mathcal{A}, \mathrm{E}[\cdot])$ be a noncommutative probability space. That is, $\mathcal{A}$ is a finite von Neumann algebra, and $\mathrm{E}[\cdot]$ is a faithful normal tracial state on it. The commutative case is included in this setting when $\mathcal{A} = L^\infty(\Lambda, P)$ for $\Lambda$ a measure space, $P$ a probability measure, and $\mathrm{E}[\cdot]$ the expectation with respect to $P$. Let $\{X(t)\}$ be an operator-valued stochastic process whose increments are stationary with respect to the state $\mathrm{E}[\cdot]$ and independent in a certain sense; see Sections 3–5 for examples of such conditions. Denote by $\mathcal{P}(n)$ the collection of all set partitions of a set of $n$ elements. For a set partition $\pi = \{B_1, B_2, \ldots, B_l\}$, temporarily denote by $c(i)$ the index of the class $B_{c(i)}$ to which $i$ belongs. Then the stochastic measure corresponding to the partition $\pi$ is

$$\mathrm{St}_\pi(t) = \int_{\substack{[0,t)^l \\ \text{all } s_i\text{'s distinct}}} dX(s_{c(1)}) \, dX(s_{c(2)}) \cdots dX(s_{c(n)}).$$

In particular, denote by $\Delta_n = \mathrm{St}_{\hat{1}}$ the higher diagonal measures of the process defined by

$$\Delta_n(t) = \int_{[0,t)} (dX(s))^n,$$

and by $\psi_n = \mathrm{St}_{\hat{0}}$, the full stochastic measures defined by

$$\psi_n(t) = \int_{\substack{[0,t)^n \\ \text{all } s_i\text{'s distinct}}} dX(s_1) \, dX(s_2) \cdots dX(s_n).$$

Here the integrals are defined by approximation, as follows.

DEFINITION 2.1. Let $\mathcal{I} = \{I_i\}_{i=1}^N$ be a subdivision of the interval $[0,t)$ into disjoint half-open intervals $I_i = [a_i, a_{i+1})$, $0 = a_1 < a_2 < \cdots < a_N < a_{N+1} = t$. Denote by $\delta(\mathcal{I}) = \max_i |I_i|$. Let $\pi \in \mathcal{P}(n)$ and $\{X(s)\}$ be a (possibly noncommutative) stochastic process. Define

$$[N]_\pi^n = \{\vec{u} \in \{1, 2, \ldots, N\}^n : u(i) = u(j) \Leftrightarrow i \overset{\pi}{\sim} j\}$$

and

$$\mathrm{St}_\pi(t; \mathcal{I}) = \sum_{\vec{u} \in [N]_\pi^n} \prod_{i=1}^n (X(a_{u(i)+1}) - X(a_{u(i)})).$$

Finally, define

$$\mathrm{St}_\pi(t) = \lim_{\delta(\mathcal{I}) \to 0} \mathrm{St}_\pi(t; \mathcal{I})$$

if the limit exists.

The existence of the limits has been proven under various conditions, see Sections 3–5 for the more precise description. For the purposes of this section we will assume that the limits exist and consider purely combinatorial facts. The most pertinent of these corresponds to linearization or, in the context of stochastic integration, to the Itô formula. Set $n = \sum_{j=1}^{k} n_j$. Denote by

$$\pi_{n_1, n_2, \ldots, n_k} \in \mathcal{P}(n)$$

the partition whose classes are intervals of consecutive integers of lengths $n_1, n_2, \ldots, n_k$. Denote

$$\mathcal{P}(n_1, n_2, \ldots, n_k) = \{\pi \in \mathcal{P}(n) : \pi \wedge \pi_{n_1, n_2, \ldots, n_k} = \hat{0}\},$$

the partitions *inhomogeneous* with respect to $\pi_{n_1, n_2, \ldots, n_k}$, that is, the collection of all partitions which do not put together elements of the $k$ distinguished subsets in the same class. For example,

$$\mathcal{P}(2,2) = \{\{(1)(2)(3)(4)\}, \{(1,3)(2)(4)\}, \{(1,4)(2)(3)\},$$
$$\{(1)(2,3)(4)\}, \{(1)(2,4)(3)\}, \{(1,3)(2,4)\}, \{(1,4)(2,3)\}\}.$$

Then

(1) $$\prod_{j=1}^{k} \psi_{n_j}(t) = \sum_{\pi \in \mathcal{P}(n_1, n_2, \ldots, n_k)} \mathrm{St}_\pi(t).$$

For a fixed subdivision $\mathcal{I}$, the statement

$$\prod_{j=1}^{k} \psi_{n_j}(t; \mathcal{I}) = \sum_{\pi \in \mathcal{P}(n_1, n_2, \ldots, n_k)} \mathrm{St}_\pi(t; \mathcal{I})$$

is purely combinatorial; see [12], Theorem 4, or [1], Proposition 4, for its proof. The fact that the relation (1) also holds in the limit will again be treated in each of the subsequent sections separately.

Denote $R_\pi(t) = \mathrm{E}[\mathrm{St}_\pi(t)]$ and $R_n(t) = \mathrm{E}[\Delta_n(t)]$. Here $R_n$ is the $n$th generalized cumulant of the process; for a process with independent increments it is the usual cumulant. Then

(2) $$\mathrm{E}\left[\prod_{j=1}^{k} \psi_{n_j}(t)\right] = \sum_{\pi \in \mathcal{P}(n_1, n_2, \ldots, n_k)} R_\pi(t).$$

For a centered process, $R_1 = 0$. In all examples we will consider, this will imply that $R_\pi = 0$ if $\pi$ contains a singleton class (a class consisting of one element). One consequence of this fact is that

(3) $$\mathrm{E}[\psi_n(t) \psi_k(t)] = 0$$

for $n \neq k$. That is, full stochastic measures of different orders are orthogonal. Thus, in general, we may consider the stochastic measures as analogs of orthogonal polynomials, and in this case formula (2) describes their linearization coefficients. The purpose of this paper is to describe examples when stochastic measures are, in fact, polynomials in the original process. If $\psi_n(t) = P_n(X(t))$, equation (3) says that the polynomials $\{P_n\}$ are orthogonal with respect to the distribution $\mu_t$ of $X(t)$ (which is a probability measure on $\mathbb{R}$). So their linearization coefficients are precisely

$$\langle P_{n_1} P_{n_2} \ldots P_{n_k} \rangle = \mathrm{E}\left[\prod_{j=1}^k P_{n_j}(X(t))\right] = \mathrm{E}\left[\prod_{j=1}^k \psi_{n_j}(t)\right].$$

Moreover, in all examples below $X(t)$ has infinite spectrum (takes on infinitely many values). So if $\mathrm{St}_\pi(t)$ is also a polynomial in $X(t)$, equation (1) implies the equality of the corresponding polynomials.

Another property which holds for some orthogonal polynomials, but which always holds for stochastic measures, is the martingale property.

PROPOSITION 2.2. *For $t > 0$, let $\mathcal{A}_t$ be the von Neumann algebra generated by the set $\{X(s) : s < t\}$. Assume the following:*

(a) *There exist consistent conditional expectations $\{\mathrm{E}_t[\cdot]\}$ from $\mathcal{A}$ onto each $\mathcal{A}_t$ preserving the expectation $\mathrm{E}[\cdot]$.*

(b) *The process $\{X(t)\}$ is centered, that is, $\mathrm{E}[X(t)] = 0$ for all $t$.*

(c) *The increments of the process are singleton independent. That is, given a collection of intervals $[s_j, t_j) \subset \mathbb{R}_+$, $j = 1, 2, \ldots, k$ such that for some $i$,*

$$[s_i, t_i) \cap \left(\bigcup_{j \neq i} [s_j, t_j)\right) = \varnothing,$$

*then $\mathrm{E}[(X(t_1) - X(s_1)) \ldots (X(t_i) - X(s_i)) \ldots (X(t_k) - X(s_k))] = 0$.*

(d) *The limit defining $\psi_n(t; X)$ exists in the $L^2$-norm with respect to $\mathrm{E}[\cdot]$.*

*Then the process $\psi_n(t; X)$ is a martingale with respect to the filtration $\{\mathcal{A}_t\}$. That is, for all $s < t$,*

$$\mathrm{E}_s[\psi(t; X)] = \psi(s; X).$$

See the Appendix for the proof.

**3. Processes with independent increments.** Let $\{X(t)\}$ be a process with stationary independent increments, and, thus, a Lévy process. Then by the results of [12], the integrals defining stochastic measures exist as limits in probability. Moreover, it is not hard to show that in this case for $\pi = \{B_1, B_2, \ldots, B_l\}$,

(4) $$\mathrm{St}_\pi(t) = \psi(t; \Delta_{|B_1|}, \Delta_{|B_2|}, \ldots, \Delta_{|B_l|}).$$

Here we are using a slightly more general definition of a stochastic measure where different factors in its defining integral may come from different processes:

$$\psi(t; (X^{(1)}, X^{(2)}, \ldots, X^{(k)})) = \int_{\substack{[0,t)^k \\ \text{all } s_i\text{'s distinct}}} dX^{(1)}(s_1) \, dX^{(2)}(s_2) \cdots dX^{(k)}(s_k).$$

See [3] for details. Throughout the paper we will consider stochastic processes for which the diagonal measures are affine functions in the original process $X$. Two types of processes that have this property are generalized Brownian motions and generalized Poisson processes.

A stochastic measure is multiplicative if $\mathrm{E}[\mathrm{St}_\pi(t)] = \prod_{B \in \pi} \mathrm{E}[\Delta_{|B|}(t)]$. Both stochastic measures in this section are multiplicative. For a multiplicative measure

$$R_\pi(t) = \prod_{B \in \pi} R_{|B|}(t), \tag{5}$$

and the sum on the right-hand side of (2) is equal to

$$\sum_{\pi \in \mathcal{P}(n_1, n_2, \ldots, n_k)} \prod_{B \in \pi} R_{|B|}(t).$$

In this case, (1) follows from [12], Theorem 4.

NOTATION 3.1. Denote by $\mathcal{P}_{1,2}(n)$ partitions whose classes consist only of one or two elements, otherwise known as "matchings," and by $\mathcal{P}_2(n)$, the collection of all pair partitions, otherwise known as "perfect matchings." Denote by $s(\pi)$ the number of singleton (one-element) classes of $\pi$, and by $s_2(\pi)$, the number of two-element classes.

3.1. *Hermite.* If $\{X(t)\}$ is the Brownian motion, then by the strong law of large numbers $\Delta_2(t) = t$, $\Delta_m(t) = 0$ for $m > 2$. Moreover, it follows from the Kailath–Segall formula (see [12], Theorem 2) that $\psi_m(t) = H_m(X(t), t)$. Here $H_m(x, t)$ are the Hermite polynomials, defined by the recursion relations

$$xH_m(x,t) = H_{m+1}(x,t) + mtH_{m-1}(x,t).$$

It follows from (4) that for $\pi \in \mathcal{P}_{1,2}$,

$$\mathrm{St}_\pi(t) = t^{s_2(\pi)} H_{s(\pi)}(X(t), t)$$

and they are 0 otherwise. Therefore, (1) gives

$$\prod_{j=1}^k H_{n_j}(x,t) = \sum_{\pi \in \mathcal{P}_{1,2}(n_1, n_2, \ldots, n_k)} t^{s_2(\pi)} H_{s(\pi)}(x,t).$$

In particular,

$$\left\langle \prod_{j=1}^{k} H_{n_j}(x,t) \right\rangle = t^{n/2}|\mathcal{P}_2(n_1,n_2,\ldots,n_k)|.$$

This formula is well known and surely quite old. Since $\langle H_m(x,t)^2\rangle = m!t^m$,

$$\prod_{j=1}^{k} H_{n_j}(x,t) = \sum_{m=0}^{n} \frac{1}{m!} t^{(n-m)/2}|\mathcal{P}_2(n_1,n_2,\ldots,n_k,m)|H_m(x,t).$$

3.2. *Centered Charlier.* If $\{X(t)\}$ is the centered Poisson process, then from [12], Proposition 7, $\Delta_m(t) = X(t) + t$ for $m \geq 2$. Moreover, $\psi_m(t) = C_m(X(t),t)$. Here $C_m(x,t)$ are the centered Charlier polynomials, defined by the recursion relations

$$xC_m(x,t) = C_{m+1}(x,t) + mC_m(x,t) + tmC_{m-1}(x,t).$$

It follows from (4) that

$$\text{St}_\pi = \sum_{l=0}^{|\pi|-s(\pi)} \binom{|\pi|-s(\pi)}{l} t^l C_{|\pi|-l}(X(t),t).$$

Therefore,

$$\prod_{j=1}^{k} C_{n_j}(x,t) = \sum_{\pi \in \mathcal{P}(n_1,n_2,\ldots,n_k)} \sum_{l=0}^{|\pi|-s(\pi)} \binom{|\pi|-s(\pi)}{l} t^l C_{|\pi|-l}(x,t).$$

In particular,

$$\left\langle \prod_{j=1}^{k} C_{n_j}(x,t) \right\rangle = \sum_{\substack{\pi \in \mathcal{P}(n_1,n_2,\ldots,n_k) \\ s(\pi)=0}} t^{|\pi|}.$$

This formula appears in [15] and a number of later sources. Since $\langle C_m(x,t)^2\rangle = m!t^m$,

$$\prod_{j=1}^{k} C_{n_j}(x,t) = \sum_{m=0}^{n} \frac{1}{m!} \sum_{\substack{\pi \in \mathcal{P}(n_1,n_2,\ldots,n_k,m) \\ s(\pi)=0}} t^{|\pi|-m} C_m(x,t).$$

Note that the noncentered polynomials, here and in the subsequent sections, will have exactly the same linearization coefficients.

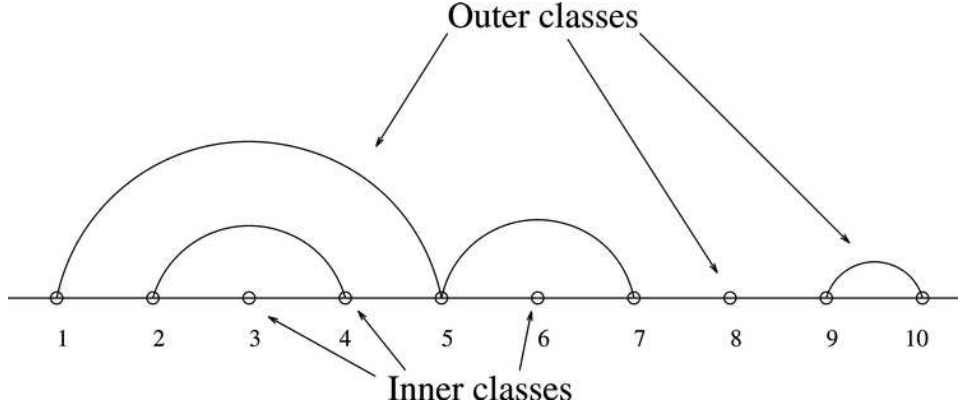

Fig. 1. *A noncrossing partition of* 10 *elements with* 3 *inner and* 3 *outer classes.*

**4. Processes with freely independent increments.** The notion of free independence was introduced by Voiculescu [14] in the context of operator algebras. $X, Y \in \mathcal{A}$ are freely independent if, whenever

$$E[f_1(X)] = E[g_1(Y)] = \cdots = E[f_n(X)] = E[g_n(Y)] = 0$$

and $g_0(Y), f_{n+1}(X)$ each are either centered or scalar, then

$$E[g_0(Y) f_1(X) g_1(Y) \ldots f_n(X) g_n(Y) f_{n+1}(X)] = 0.$$

This property is easily seen to be incompatible with, but is parallel to, the usual independence. Free probability is by now quite a rich theory which is based on this notion; see [14] for an overview. In particular, there is a well-developed theory of free cumulants, free infinitely divisible distributions and limit theorems, and processes with freely independent increments.

In this section, let $\{X(t)\}$ be a process with stationary freely independent increments, and, thus, a free Lévy process. It was shown in [1, 3] that in this case the integrals defining stochastic measures exist as limits in the operator norm, as long as the operators $\{X(t)\}$ are bounded. Moreover, it was shown in [1] that $\mathrm{St}_\pi = 0$ unless $\pi$ is a noncrossing partition. Here a partition $\pi$ is noncrossing if there are no $i < j < k < l$ with $i \overset{\pi}{\sim} k$, $j \overset{\pi}{\sim} l$, $i \overset{\pi}{\not\sim} j$.

In this case, (1) follows from [1], Proposition 4. For the analog of the formula (4), we need a new notion. For a noncrossing partition, we distinguish the classes that are inner, or covered by other classes, and outer. See Figure 1 for an example.

PROPOSITION 4.1. *Let $\pi$ be a noncrossing partition with outer classes $B_1, \ldots, B_{o(\pi)}$ and inner classes $C_1, \ldots, C_{i(\pi)}$. Then*

$$\mathrm{St}_\pi(t) = \prod_{i=1}^{i(\pi)} R_{|C_i|}(t) \cdot \psi(\Delta_{|B_1|}(t), \Delta_{|B_2|}(t), \ldots, \Delta_{|B_{o(\pi)}|}(t)).$$

PROOF. This is a particular case of the main theorem of [3]. □

NOTATION 4.2. Denote by $NC(n)$ the lattice of noncrossing partitions, and by $NC_{1,2}(n)$, $NC_2(n)$, $NC(n_1,\ldots,n_k)$, and so on, the corresponding subsets of $NC(n)$. Denote by $si(\pi)$ the number of inner singletons of $\pi$. Denote by $o(\pi)$ and $i(\pi)$ the number of outer and, respectively, inner classes of $\pi$.

Free stochastic measures are not multiplicative in general. However, (5) does hold for $\pi \in NC(n)$. So for a free stochastic measure, the general linearization coefficients are

$$\mathrm{E}\left[\prod_{j=1}^{k}\psi_{n_j}(t)\right] = \sum_{\pi\in NC(n_1,n_2,\ldots,n_k)} R_\pi(t) = \sum_{\pi\in NC(n_1,n_2,\ldots,n_k)} \prod_{B\in\pi} R_{|B|}(t).$$

4.1. *Chebyshev.* There is a free version of the central limit theorem, with independent variables replaced by freely independent ones. The limit distribution in this theorem is the semicircular distribution. A process (consisting of noncommuting operators) $\{X(t)\}$ with stationary freely independent increments all of which have (scaled) semicircular distributions is the free Brownian motion.

If $\{X(t)\}$ is the free Brownian motion, then from [1], $\Delta_2(t) = t$, $\Delta_m(t) = 0$ for $m > 2$. Moreover, by [1], Corollary 8, $\psi_m(t) = U_m(X(t),t)$. Here $U_m(x,t)$ are the Chebyshev polynomials of the second kind, defined by the recursion relations

$$xU_m(x,t) = U_{m+1}(x,t) + tU_{m-1}(x,t).$$

It follows from Proposition 4.1 that for $\pi \in NC_{1,2}(n)$ and $si(\pi) = 0$,

$$\mathrm{St}_\pi(t) = t^{s_2(\pi)}U_{s(\pi)}(X(t),t)$$

and they are 0 otherwise. Therefore, by (1),

$$\prod_{j=1}^{k}U_{n_j}(x,t) = \sum_{\substack{\pi\in NC_{1,2}(n_1,n_2,\ldots,n_k)\\si(\pi)=0}} t^{s_2(\pi)}U_{s(\pi)}(x,t).$$

In particular,

$$\left\langle\prod_{j=1}^{k}U_{n_j}(x,t)\right\rangle = t^{n/2}|NC_2(n_1,n_2,\ldots,n_k)|.$$

This formula has essentially appeared in [8], in a slightly different guise (they count the number of Dyck paths). Since $\langle U_m(x,t)^2\rangle = t^m$,

$$\prod_{j=1}^{k}U_{n_j}(x,t) = \sum_{m=0}^{n} t^{(n-m)/2}|NC_2(n_1,n_2,\ldots,n_k,m)|U_m(x,t).$$

4.2. *Centered free Charlier.* The distribution of the sum on $n$ freely independent Bernoulli $((1 - \frac{t}{n})\delta_0 + \frac{t}{n}\delta_1)$ variables converges, as $n \to \infty$, to a distribution which is naturally called the free Poisson distribution with parameter $t$. It is also known as the Marchenko–Pastur (or, for $t = 1$, Wishart) distribution. A process with stationary freely independent increments such that the increments have free Poisson distributions is the free Poisson process.

If $\{X(t)\}$ is the centered free Poisson process, then by [1], Corollary 4, $\Delta_m(t) = X(t) + t$ for $m \geq 2$. By [1], Corollary 10, $\psi_m(t) = C_{0,m}(X(t), t)$. Here $C_{0,m}(x, t)$ are the centered free Charlier polynomials, defined by the recursion relations

$$xC_{0,0}(x,t) = C_{0,1}(x,t),$$
$$xC_{0,m}(x,t) = C_{0,m+1}(x,t) + C_{0,m}(x,t) + tC_{0,m-1}(x,t)$$

for $m > 0$. They are, of course, orthogonal with respect to the free Poisson distribution.

It follows from Proposition 4.1 that for $\pi \in NC(n)$ and $si(\pi) = 0$,

$$\mathrm{St}_\pi = t^{i(\pi)} \sum_{l=0}^{o(\pi)-s(\pi)+si(\pi)} \binom{o(\pi)-s(\pi)+si(\pi)}{l} t^l C_{0,o(\pi)-l}(X(t),t)$$

and they are 0 otherwise. Therefore, by (1),

$$\prod_{j=1}^{k} C_{0,n_j}(x,t)$$

$$= \sum_{\substack{\pi \in NC(n_1,n_2,\ldots,n_k),\\ si(\pi)=0}} \sum_{l=0}^{o(\pi)-s(\pi)+si(\pi)} \binom{o(\pi)-s(\pi)+si(\pi)}{l}$$

$$\times t^{i(\pi)+l} C_{0,o(\pi)-l}(x,t).$$

In particular,

$$\left\langle \prod_{j=1}^{k} C_{0,n_j}(x,t) \right\rangle = \sum_{\substack{\pi \in NC(n_1,n_2,\ldots,n_k)\\ s(\pi)=0}} t^{|\pi|}.$$

Since $\langle C_{0,m}(x,t)^2 \rangle = t^m$,

$$\prod_{j=1}^{k} C_{0,n_j}(x,t) = \sum_{m=0}^{n} \sum_{\substack{\pi \in NC(n_1,n_2,\ldots,n_k,m)\\ s(\pi)=0}} t^{|\pi|-m} C_{0,m}(x,t).$$

## 5. Processes on a $q$-deformed full Fock space.

### 5.1. $q$-Fock space.
Consider the Hilbert space $L^2(\mathbb{R}_+, dx)$. Let

$$\mathcal{F}_{\text{alg}}(L^2(\mathbb{R}_+)) = \bigoplus_{k=0}^{\infty} L^2(\mathbb{R}_+, dx)^{\otimes k} = \bigoplus_{k=0}^{\infty} L^2(\mathbb{R}_+^k, dx^{\otimes k})$$

be its algebraic Fock space. Here the 0th component is spanned by the vacuum vector $\Omega$. Then $\langle \cdot, \cdot \rangle_0$ defined by

$$\langle f_1 \otimes \cdots \otimes f_k, g_1 \otimes \cdots \otimes g_n \rangle_0 = \delta_{kn} \langle f_1, g_1 \rangle \cdots \langle f_k, g_k \rangle$$

is an inner product on the algebraic Fock space, where $\langle \cdot, \cdot \rangle$ is the standard inner product on $L^2(\mathbb{R}_+, dx)$. Define the operator $P_q$ by

$$P_q(f_1 \otimes \cdots \otimes f_n) = \sum_{\sigma \in \text{Sym}(n)} q^{i(\sigma)} f_{\sigma(1)} \otimes \cdots \otimes f_{\sigma(n)},$$

where $\text{Sym}(n)$ is the permutation group and $i(\sigma)$ is the number of inversions of $\sigma$. According to [7], this operator is strictly positive for $-1 < q < 1$. Denote $\langle \cdot, \cdot \rangle_q = \langle \cdot, P_q \cdot \rangle_0$. Then this is also an inner product, and we denote by $\mathcal{F}_q(L^2(\mathbb{R}_+))$ the completion of $\mathcal{F}_{\text{alg}}(L^2(\mathbb{R}_+))$ with respect to the corresponding norm, and call it the $q$-deformed full Fock space.

For $f \in L^2(\mathbb{R}_+) \cap L^{\infty}(\mathbb{R}_+)$, define creation, annihilation and preservation operators on the $q$-Fock space $\mathcal{F}_q(L^2(\mathbb{R}_+))$ by

$$a^*(f)(\Omega) = f,$$

$$a^*(f)(g_1 \otimes \cdots \otimes g_n) = f \otimes g_1 \otimes \cdots \otimes g_n,$$

$$a(f)(\Omega) = 0,$$

$$a(f)(g) = \langle f, g \rangle \Omega,$$

$$a(f)(g_1 \otimes \cdots \otimes g_n) = \sum_{k=1}^{n} q^{k-1} \langle f, g_k \rangle g_1 \otimes \cdots \otimes \hat{g}_k \otimes \cdots \otimes g_n,$$

$$p(f)(\Omega) = 0,$$

$$p(f)(g_1 \otimes \cdots \otimes g_n) = \sum_{k=1}^{n} q^{k-1} f g_k \otimes g_1 \otimes \cdots \otimes \hat{g}_k \otimes \cdots \otimes g_n,$$

where $\hat{g}_k$ means "omit $k$th term." For $f$ real-valued, $p(f)$ is self-adjoint, and $a(f)$ and $a^*(f)$ are adjoints of each other.

The noncommutative stochastic process

$$X(t) = a^*(\mathbf{1}_{[0,t)}) + a(\mathbf{1}_{[0,t)})$$

is, by definition, the $q$-Brownian motion, and the process

$$X(t) = a^*(\mathbf{1}_{[0,t)}) + a(\mathbf{1}_{[0,t)}) + p(\mathbf{1}_{[0,t)})$$

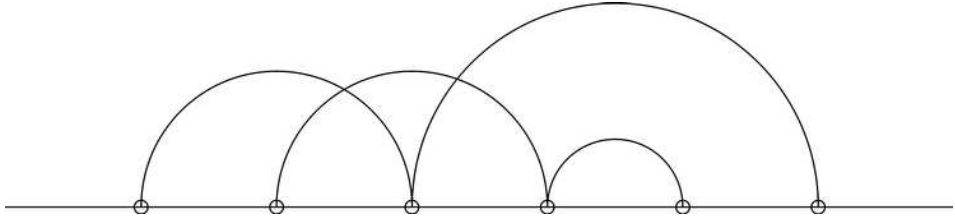

Fig. 2. *A partition of 6 elements with 2 restricted crossings.*

is the centered $q$-Poisson process. Let $\mathcal{A}$ be the von Neumann algebra generated by $\{X(t)\}_{t \in [0,\infty)}$, and let $\mathrm{E}[\cdot]$ be the vacuum vector state $\mathrm{E}[\cdot] = \langle \Omega, \cdot \Omega \rangle$. Then for $q = 0$, these processes are the free Brownian motion and the centered free Poisson process, while for the degenerate case $q = 1$, they give the corresponding classical processes.

NOTATION 5.1. Let $\pi$ be a partition. Denote by $\mathrm{rc}(\pi)$ the number of restricted crossings of $\pi$. Here a restricted crossing of $\pi$ is a 4-tuple $i < j < k < l$ such that $i \stackrel{\pi}{\sim} k$, $j \stackrel{\pi}{\sim} l$, and $k = \min\{r : r > i, r \stackrel{\pi}{\sim} i\}$, $l = \min\{r : r > j, r \stackrel{\pi}{\sim} j\}$. See Figure 2 for an example. Also, define the singleton depth $\mathrm{sd}(\pi)$ of $\pi$ to be the sum of depths, $d(i) = |\{j | \exists a, b \in B_j : a < i < b\}|$, over all the singletons $(i)$ of $\pi$.

For $-1 < q < 1$, denote

$$[0]_q = 0, \qquad [n]_q = \sum_{j=0}^{n-1} q^j = \frac{1-q^n}{1-q} \quad \text{and} \quad [n]_q! = \prod_{j=1}^{n} [j]_q.$$

The stochastic measures for the $q$-Lévy processes are described in a forthcoming paper [4]. However, the functionals $R_\pi$ are known to be well defined, and the following analog of (5) holds.

PROPOSITION 5.2 ([2], Theorem 3.8).

$$R_\pi(t) = q^{\mathrm{rc}(\pi)} \prod_{B \in \pi} R_{|B|}(t).$$

Therefore, if (1) holds, the linearization coefficients are

$$\mathrm{E}\left[\prod_{j=1}^{k} \psi_{n_j}(t)\right] = \sum_{\pi \in \mathcal{P}(n_1, n_2, \ldots, n_k)} R_\pi(t) = \sum_{\pi \in \mathcal{P}(n_1, n_2, \ldots, n_k)} q^{\mathrm{rc}(\pi)} \prod_{B \in \pi} R_{|B|}(t).$$

PROPOSITION 5.3. *If $X(t)$ is the $q$-Brownian motion, then the limit defining $\Delta_k(t; X)$ exists in the $L^2$-norm with respect to $\mathrm{E}[\cdot]$, and equals $\Delta_2(t) = t$, $\Delta_k(t) = 0$ for $k > 2$. Similarly, if $X(t)$ is the centered $q$-Poisson* process, the limit defining $\Delta_k(t; X)$ exists in the $L^2$-norm with respect to $\mathrm{E}[\cdot]$, and equals $\Delta_k(t) = X(t) + t$ for $k \geq 2$.

See the Appendix for the proof.

5.2. *q-Hermite.* If $\{X(t)\}$ is the $q$-Brownian motion, then $\psi_m(t) = H_{q,m}(X(t), t)$. Here the polynomials $H_{q,m}(x, t)$ are a scaled version of the continuous (Rogers) $q$-Hermite polynomials, defined by the recursion relations

$$xH_{q,m}(x,t) = H_{q,m+1}(x,t) + t[m]_q H_{q,m-1}(x,t).$$

The measure of orthogonality of these polynomials, and so the distribution of $X(t)$, is the most common version of the $q$-Gaussian measure; see [6] or (with a slightly different normalization) [11], 3.26.

LEMMA 5.4 ([2], Proposition 6.12). *For $\pi \in \mathcal{P}_{1,2}$,*

$$\mathrm{St}_\pi = q^{\mathrm{rc}(\pi)+\mathrm{sd}(\pi)} t^{s_2(\pi)} H_{q,s(\pi)}(X(t), t)$$

*and they are 0 otherwise. Here the limit in the definition of the stochastic measure is in $L^{\infty-}$, that is, in $L^p$ for any $1 \leq p < \infty$, with respect to $\mathrm{E}[\cdot]$.*

Therefore,

$$\prod_{j=1}^k H_{q,n_j}(x,t) = \sum_{\pi \in \mathcal{P}_{1,2}(n_1,n_2,\ldots,n_k)} q^{\mathrm{rc}(\pi)+\mathrm{sd}(\pi)} t^{s_2(\pi)} H_{q,s(\pi)}(x,t).$$

In particular,

$$\left\langle \prod_{j=1}^k H_{q,n_j}(x,t) \right\rangle = t^{n/2} \sum_{\pi \in \mathcal{P}_2(n_1,n_2,\ldots,n_k)} q^{\mathrm{rc}(\pi)}.$$

This formula has appeared in [9]. Since $\langle H_{q,m}(x,t)^2 \rangle = [m]_q! t^m$,

$$\prod_{j=1}^k H_{q,n_j}(x,t) = \sum_{m=0}^n \frac{1}{[m]_q!} t^{(n-m)/2} \sum_{\pi \in \mathcal{P}_2(n_1,n_2,\ldots,n_k,m)} q^{\mathrm{rc}(\pi)} H_{q,m}(x,t).$$

5.3. *Centered big q-Hermite.* Let $\{X(t)\}$ be the centered $q$-Poisson process.

PROPOSITION 5.5. *For the centered $q$-Poisson process, $\psi_m(t) = C_{q,m}(X(t), t)$. Here $C_{q,m}$ are a scaled version of the centered continuous big $q$-Hermite polynomials, which in our context are $q$-analogs of the Charlier polynomials. They are defined by the recursion relations*

$$xC_{q,m}(x,t) = C_{q,m+1}(x,t) + [m]_q C_{q,m}(x,t) + t[m]_q C_{q,m-1}(x,t).$$

*In particular, the stochastic measures $\psi_m$ are well defined, with the limits taken in the $L^2$-norm with respect to $\mathrm{E}[\cdot]$.*

See the Appendix for the proof. Also, see [13] for a detailed description of the measure of orthogonality of these polynomials.

From Proposition 5.2, for $s(\pi) = 0$,

$$R_\pi(t) = q^{\mathrm{rc}(\pi)} t^{|\pi|}, \tag{6}$$

and they are 0 otherwise.

PROOF OF (1). We start with the known combinatorial formula

$$\prod_{j=1}^{k} \psi_{n_j}(t; \mathcal{I}) = \sum_{\pi \in \mathcal{P}(n_1, n_2, \ldots, n_k)} \mathrm{St}_\pi(t; \mathcal{I}).$$

By the results of [4], $\mathrm{St}_\pi(t; \mathcal{I})$ converges to $\mathrm{St}_\pi(t)$ in $L^2$. By the previous proposition, $\psi_{n_j}(t; \mathcal{I})$ converges to $\psi_{n_j}(t)$ in $L^2$. Thus, $\psi_{n_1}(t; \mathcal{I})\psi_{n_2}(t; \mathcal{I})$ converges to $\psi_{n_1}(t)\psi_{n_2}(t)$ in $L^1$. On the other hand, it also converges to $\sum_{\pi \in \mathcal{P}(n_1, n_2)} \mathrm{St}_\pi(t)$ in $L^2$. Therefore, the two latter expressions are equal, and $\psi_{n_1}(t; \mathcal{I})\psi_{n_2}(t; \mathcal{I})$ converges, in fact, in $L^2$. By induction, we conclude that $\prod_{j=1}^{k} \psi_{n_j}(t; \mathcal{I})$ converges in $L^2$ to $\sum_{\pi \in \mathcal{P}(n_1, n_2, \ldots, n_k)} \mathrm{St}_\pi(t)$. □

We conclude that

$$\left\langle \prod_{j=1}^{k} C_{q,n_j}(x,t) \right\rangle = \sum_{\substack{\pi \in \mathcal{P}(n_1, n_2, \ldots, n_k) \\ s(\pi) = 0}} q^{\mathrm{rc}(\pi)} t^{|\pi|}.$$

Since $\langle C_{q,m}(x,t)^2 \rangle = [m]_q! t^m$,

$$\prod_{j=1}^{k} C_{q,n_j}(x,t) = \sum_{m=0}^{n} \frac{1}{[m]_q!} \sum_{\substack{\pi \in \mathcal{P}(n_1, n_2, \ldots, n_k, m) \\ s(\pi) = 0}} q^{\mathrm{rc}(\pi)} t^{|\pi|-m} C_{0,m}(x,t).$$

5.4. *Limiting relations.* The results of the previous sections can be obtained as the limits of the results of this one. For the continuous (Rogers) $q$-Hermite polynomials, taking $q = 1$ gives the formulas for the Hermite polynomials, while taking $q = 0$ gives the formulas for the Chebyshev polynomials. For the continuous big $q$-Hermite polynomials, taking $q = 1$ gives the formulas for the Charlier polynomials, while taking $q = 0$ gives the formulas for the free Charlier polynomials. Note that in the latter case we only recover the linearization coefficients themselves, not the expressions for the products of polynomials as sums over partitions.

Finally, consider the process

$$X(t, \alpha) = a^*(\mathbf{1}_{[0,t)}) + a(\mathbf{1}_{[0,t)}) + \alpha p(\mathbf{1}_{[0,t)}).$$

For this process, $\Delta_m(t) = \alpha^{m-1} X(t) + \alpha^{m-2} t$ and $R_m(t) = \alpha^{m-2} t$ for $m \geq 2$. Therefore, for $s(\pi) = 0$,

$$R_\pi = q^{\mathrm{rc}(\pi)} t^{|\pi|} \alpha^{n-2|\pi|},$$

and they are 0 otherwise. Also, for this process, $\psi_m(t) = P_{q,m,\alpha}(X(t), t)$, where

$$x P_{q,m,\alpha}(x,t) = P_{q,m+1,\alpha}(x,t) + \alpha [m]_q P_{q,m,\alpha}(x,t) + t[m]_q P_{q,m-1,\alpha}(x,t).$$

We conclude that

$$\left\langle \prod_{j=1}^k P_{q,n_j,\alpha}(x,t) \right\rangle = \sum_{\substack{\pi \in \mathcal{P}(n_1, n_2, \ldots, n_k) \\ s(\pi) = 0}} q^{\mathrm{rc}(\pi)} t^{|\pi|} \alpha^{n-2|\pi|}.$$

For $\alpha = 1$, this gives the $q$-Poisson process and the continuous big $q$-Hermite polynomials. On the other hand, for $\alpha = 0$, this gives the $q$-Brownian motion and the continuous (Rogers) $q$-Hermite polynomials. In the linearization formula, the only partitions with a nonzero contribution are those with $n = 2|\pi|$ and without singletons, that is, pair partitions.

## APPENDIX

**$q$-Lévy processes.** We briefly review the definition of more general $q$-Lévy processes and their stochastic measures; see [2] for more details. Let $V$ be a Hilbert space, and consider $H = L^2(\mathbb{R}_+, dx) \otimes V$. Define $\mathcal{F}_{\mathrm{alg}}(H)$, $\mathcal{F}_q(H)$, E$[\cdot]$ and, for $\xi \in H$, $a(\xi)$ and $a^*(\xi)$ as in the beginning of Section 5 for $V = \mathbb{C}$. For $T$ an essentially self-adjoint operator on $H$, define $p(T)$ on $\mathcal{F}_q(H)$ by

$$p(T)(\Omega) = 0,$$

$$p(T)(\xi_1 \otimes \cdots \otimes \xi_n) = \sum_{k=1}^n q^{k-1} (T\xi_k) \otimes \xi_1 \otimes \cdots \otimes \hat{\xi}_k \otimes \cdots \otimes \xi_n.$$

By [2], Proposition 2.2, $p(T)$ is an essentially self-adjoint operator.

Pick $\xi \in V$, $T$ an operator on $V$ and $\lambda \in \mathbb{R}$. Assume that

(1)      $T$ is essentially self-adjoint, the vectors $\{T^k \xi\}_{k=0}^\infty$ belong to its dense domain, span it, and are analytic for $T$.

Define $a_t(\xi) = a(\mathbf{1}_{[0,t)} \otimes \xi)$, $a_t^*(\xi) = a^*(\mathbf{1}_{[0,t)} \otimes \xi)$, and $p_t(T) = p(\mathbf{1}_{[0,t)} \otimes T)$. Then the corresponding $q$-Lévy process is

$$p_t(\xi, T, \lambda) = a_t^*(\xi) + a_t(\xi) + p_t(T) + \lambda t.$$

Let $\{X(t)\}$ be such a process, let $\mathcal{I}$ be a subdivision of the interval $[0,t)$, and let $\{X_i\}$ be the increments of this process corresponding to the subdivision intervals of $\mathcal{I}$, $X_i = X(a_{i+1}) - X(a_i)$ for $I_i = [a_i, a_{i+1}) \in \mathcal{I}$. Then

$$\Delta_k(t; X, \mathcal{I}) = \sum_i X_i^k$$

and

$$\psi_k(t; X, \mathcal{I}) = \sum_{\substack{i_1, i_2, \ldots, i_k \\ \text{distinct}}} X_{i_1} \ldots X_{i_k}.$$

The stochastic measures $\Delta_k(t; X)$ and $\psi_k(t; X)$ are the limits of the above quantities as the size of the subdivision $\delta(\mathcal{I})$ tends to 0, if these limits exist.

Similarly, for a $k$-tuple of processes $(X^{(1)}, X^{(2)}, \ldots, X^{(k)})$, we can define

$$\Delta(t; (X^{(1)}, X^{(2)}, \ldots, X^{(k)}), \mathcal{I}) = \sum_i X_i^{(1)} X_i^{(2)} \ldots X_i^{(k)}.$$

Such a $k$-tuple form a multiple $q$-Lévy process if they satisfy an extra compatibility condition ([2], equation 1), similar to the one in (1).

LEMMA A.1 ([2], Proposition 3.6).  *Let $\{X^{(i)}(t) = p_t(\xi_i, T_i, \lambda_i)\}_{i=1}^k$ be a multiple $q$-Lévy process. Then the $q$-cumulants*

$$R(t; (X^{(1)}, X^{(2)}, \ldots, X^{(k)})) = \lim_{\delta(\mathcal{I}) \to 0} \mathrm{E}[\Delta(t; (X^{(1)}, X^{(2)}, \ldots, X^{(k)}), \mathcal{I})]$$

*are well defined, and equal to*

$$R(t; (X^{(1)}, X^{(2)}, \ldots, X^{(k)})) = \begin{cases} t\lambda_1, & \text{if } k = 1, \\ t\left\langle \xi_1, \prod_{j=2}^{k-1} T_j \xi_k \right\rangle, & \text{if } k \geq 2. \end{cases}$$

PROPOSITION A.2.  *Let $X(t) = p_t(\xi, T, \lambda)$ be a general one-dimensional $q$-Lévy process. Then the limit defining $\Delta_k(t; X)$ exists in the $L^2$-norm with respect to $\mathrm{E}[\cdot]$, and equals*

$$Y(t) = p_t(T^{k-1}\xi, T^k, \langle \xi, T^{k-2}\xi \rangle).$$

PROOF.  Condition (1) implies that any $k$-tuple of processes whose components are $X$ and $Y$ is also compatible. It suffices to show that

$$\lim_{\delta(\mathcal{I}) \to 0} \langle (\Delta_k(t; X, \mathcal{I}) - Y(t))^2 \Omega, \Omega \rangle = 0.$$

First expand

$$\left(\sum_i X_i^k - Y(t)\right)^2$$

$$= \sum_i X_i^{2k} + \sum_{i\neq j} X_i^k X_j^k - \sum_i Y_i X_i^k - \sum_{i\neq j} Y_i X_j^k - \sum_i X_i^k Y_i - \sum_{i\neq j} X_i^k Y_j + Y^2$$

$$= \Delta_{2k}(t;X,\mathcal{I}) + \sum_{i\neq j} X_i^k X_j^k - \Delta(t;(Y,X,\ldots,X),\mathcal{I}) - \sum_{i\neq j} Y_i X_j^k$$

$$- \Delta(t;(X,\ldots,X,Y),\mathcal{I}) - \sum_{i\neq j} X_i^k Y_j + Y(t)^2.$$

From the pyramidal independence of the increments ([2], Lemma 3.3), it follows that

$$(2) \quad \mathrm{E}\left[\sum_{i\neq j} X_i^k X_j^k\right] = \sum_{i\neq j} \mathrm{E}[X_i^k]\mathrm{E}[X_j^k] = R_k(t;X,\mathcal{I})^2 - \sum_i \mathrm{E}[X_i^k]^2.$$

By combining Proposition 5.2 with Lemma A.1 and the moment-cumulant formula ([2], equation 3)

$$\mathrm{E}[X(t)^k] = \sum_{\pi \in \mathcal{P}(k)} R_\pi(t),$$

it follows that $\mathrm{E}[X(t)^k] = \sum_{\pi \in \mathcal{P}(k)} t^{|\pi|} R_\pi(1) = O(t)$. Since the increments of the process are stationary, we also know that $\mathrm{E}[X_i^k] = \mathrm{E}[X(|I_i|)^k]$. Therefore,

$$\lim_{\delta(\mathcal{I})\to 0} \sum_i \mathrm{E}[X_i^k]^2 \leq C \lim_{\delta(\mathcal{I})\to 0} \sum_i |I_i|^2 \leq C \lim_{\delta(\mathcal{I})\to 0} t\delta(\mathcal{I}) = 0.$$

We conclude that the limit of the expression (2) is $R_k(t;X)^2$. Similarly,

$$\lim_{\delta(\mathcal{I})\to 0} \mathrm{E}\left[\sum_{i\neq j} Y_i X_j^k\right] = \lim_{\delta(\mathcal{I})\to 0} \mathrm{E}\left[\sum_{i\neq j} X_j^k Y_i\right] = \mathrm{E}[Y(t)]R_k(t;X).$$

Therefore,

$$\lim_{\delta(\mathcal{I})\to 0} \langle (\Delta_k(t;X,\mathcal{I}) - Y(t))^2 \Omega, \Omega \rangle$$

$$= R_{2k}(t;X) + R_k(t;X)^2 - R(t;(Y,X,\ldots,X)) - \mathrm{E}[Y(t)]R_k(t;X)$$
$$\quad - R(t;(X,\ldots,X,Y)) - R_k(t;X)\mathrm{E}[Y(t)] + \mathrm{E}[Y(t)^2]$$
$$= \langle \xi, T^{2k-2}\xi \rangle + \langle \xi, T^{k-2}\xi \rangle^2 - \langle T^{k-1}\xi, T^{k-1}\xi \rangle - \langle \xi, T^{k-2}\xi \rangle\langle \xi, T^{k-2}\xi \rangle$$
$$\quad - \langle \xi, T^{k-1}T^{k-1}\xi \rangle - \langle \xi, T^{k-2}\xi \rangle\langle \xi, T^{k-2}\xi \rangle$$
$$\quad + \langle T^{k-1}\xi, T^{k-1}\xi \rangle + \langle \xi, T^{k-2}\xi \rangle^2$$
$$= 0.$$

□

PROOF OF PROPOSITION 5.3. For the $q$-Brownian motion, $T = 0$, while for the $q$-Poisson process, $T = \text{Id}$. So the result follows from the preceding proposition. □

PROOF OF PROPOSITION 5.5. It suffices to show that the stochastic measures satisfy the same recursion relations as the orthogonal polynomials $C_{q,n}$. That is, we will show that

$$\text{(3)} \quad \lim_{\delta(\mathcal{I}) \to 0} \| X(t)\psi_n(t; X, \mathcal{I}) - \psi_{n+1}(t; X, \mathcal{I}) \\ - [n]_q \psi_n(t; X, \mathcal{I}) - t[n]_q \psi_{n-1}(t; X, \mathcal{I}) \|_2 = 0.$$

Indeed, if that is the case, then

$$\psi_{n+1}(t; X, \mathcal{I}) = X(t)\psi_n(t; X, \mathcal{I}) + [n]_q \psi_n(t; X, \mathcal{I}) + t[n]_q \psi_{n-1}(t; X, \mathcal{I}) + A(\mathcal{I}),$$

with $L^2 - \lim_{\delta(\mathcal{I}) \to 0} A(\mathcal{I}) = 0$. By induction, the right-hand side converges in $L^2$ to

$$X(t)\psi_n(t; X) + [n]_q \psi_n(t; X) + t[n]_q \psi_{n-1}(t; X)$$

[for the first term, we use an argument similar to the proof of (1)]. So the left-hand side also converges in $L^2$. Moreover, also by induction, the limit of the right-hand side is

$$X(t)C_{q,n}(X(t)) + [n]_q C_{q,n}(X(t)) + t[n]_q C_{q,n-1}(X(t)) = C_{q,n+1}(X(t)).$$

We will omit $X, \mathcal{I}$ and $t$ in the notation. Expanding the norm in (3), we get

$$\text{(4)} \quad \begin{aligned} &\text{E}[(\psi_{n+1} + [n]_q \psi_n + [n]_q t\psi_{n-1} - X\psi_n) \\ &\quad \times (\psi_{n+1} + [n]_q \psi_n + [n]_q t\psi_{n-1} - \psi_n X)] \\ &= \text{E}[\psi_{n+1}\psi_{n+1}] + [n]_q^2 \text{E}[\psi_n \psi_n] \\ &\quad + [n]_q^2 t^2 \text{E}[\psi_{n-1}\psi_{n-1}] + \text{E}[X\psi_n \psi_n X] \\ &\quad - \text{E}[\psi_{n+1}\psi_n X] - \text{E}[X\psi_n \psi_{n+1}] - [n]_q \text{E}[\psi_n \psi_n X] \\ &\quad - [n]_q \text{E}[X\psi_n \psi_n] - [n]_q t \text{E}[\psi_{n-1}\psi_n X] - [n]_q t \text{E}[X\psi_n \psi_{n-1}]. \end{aligned}$$

Combining the general linearization formula (2) and the specific form (6) of the cumulants of the $q$-Poisson process,

$$\text{E}[\psi_n \psi_n] = \sum_{\pi \in \mathcal{P}_2(n,n)} q^{\text{rc}(\pi)} t^n.$$

It is easy to show by induction that

$$\sum_{\pi \in \mathcal{P}_2(n,n)} q^{\text{rc}(\pi)} = [n]_q!.$$

This follows from the fact that 1 is connected by $\pi$ to exactly one element $2n-k$, their class crosses exactly $k-1$ other classes, and $\sum_{k=1}^{n} q^{k-1} = [n]_q$. We conclude that

$$\mathrm{E}[\psi_n \psi_n] = [n]_q! t^n.$$

We similarly simplify the other expressions in the sum (4). We treat in detail the most complicated term

$$\mathrm{E}[X\psi_n\psi_n X] = t\mathrm{E}[\psi_n\psi_n] + [n]_q^2 \mathrm{E}[\psi_n\psi_n] + (1+q)[n]_q^2 t^2 \mathrm{E}[\psi_{n-1}\psi_{n-1}]$$
$$= [n]_q! t^{n+1} + [n]_q^2 [n]_q! t^n + (1+q)[n]_q[n]_q! t^{n+1}.$$

Denote by $\vee$ the join of set partitions. In the first three terms in the sum below, $\pi$ is a partition in $\mathcal{P}_2(n,n)$ induced on the subset $\{2,\ldots,2n+1\}$ of $\{1,\ldots,2n+2\}$. In the last two terms, $\pi$ is a partition in $\mathcal{P}_2(n-1,n-1)$ induced on the subset $\{2,\ldots,\widehat{1+k},\ldots,\widehat{2n+2-j},\ldots,2n+1\}$ of $\{1,\ldots,2n+2\}$. Using (1),

$$\psi_1 \psi_n \psi_n \psi_1$$
$$= \sum_{\substack{\pi \in \mathcal{P}_2(n,n) \\ \pi \upharpoonright \{2,\ldots,2n+1\}}} \mathrm{St}_{\pi \vee \{(1,2n+2),(2),\ldots,(2n+1)\}}$$
$$+ \sum_{\substack{\pi \in \mathcal{P}_2(n,n) \\ \pi \upharpoonright \{2,\ldots,2n+1\}}} \sum_{\substack{k,j=1,\ldots,n \\ k \sim j \mod \pi}} \mathrm{St}_{\pi \vee \{(1,1+k),(2n+2-j,2n+2),(2),\ldots,(2n+1)\}}$$
$$+ \sum_{\substack{\pi \in \mathcal{P}_2(n,n) \\ \pi \upharpoonright \{2,\ldots,2n+1\}}} \sum_{\substack{k,j=1,\ldots,n \\ k \not\sim j \mod \pi}} \mathrm{St}_{\pi \vee \{(1,1+k),(2n+2-j,2n+2),(2),\ldots,(2n+1)\}}$$
$$+ \sum_{k,j=1}^{n} \sum_{\substack{\pi \in \mathcal{P}(n-1,n-1) \\ \pi \upharpoonright \{2,\ldots,\widehat{1+k},\ldots,\widehat{2n+2-j},\ldots,2n+1\}}} \mathrm{St}_{\pi \vee \{(1,1+k),(2n+2-j,2n+2),(2),\ldots,(2n+1)\}}$$
$$+ \sum_{k,j=1}^{n} \sum_{\substack{\pi \in \mathcal{P}(n-1,n-1) \\ \pi \upharpoonright \{2,\ldots,\widehat{1+k},\ldots,\widehat{2n+2-j},\ldots,2n+1\}}} \mathrm{St}_{\pi \vee \{(1,2n+2-j),(1+k,2n+2),(2),\ldots,(2n+1)\}}$$

$+$ terms containing singletons.

See Figure 3 for an illustration. Five types of partitions in it correspond to the five terms in the expression above. The classes containing the first (i.e., "1") and the last (i.e., "$2n+2$") elements of the set are shown; the remaining classes consist of two elements each and are inhomogeneous with respect to the partition $\pi_{n+1,n+1}$. The dashed classes belong to the partition $\pi$. Note that in the third diagram, the dashed lines may cross each other as shown in the diagram, in which case the crossing is counted among the crossings of $\pi$. Alternatively, the dashed lines may not cross, in which case one of them crosses a solid line, and this crossing is counted among the extra ones in the sum below. In the first diagram, there are no crossings between classes of $\pi$ and the extra class. In the next three diagrams, the two extra classes cover $(k-1)$ and $(j-1)$ points, respectively. Because $\pi$ is inhomogeneous, each point covered by an extra class has to be connected to a point not covered by it, hence, extra $(k-1)+(j-1)$ crossings are introduced. In the last diagram, the classes containing points $\{2,\ldots,k\}$ have to cross the extra class $(1+k, 2n+2)$, and the classes containing points $\{2n+3-j,\ldots,2n+1\}$ have to cross the extra class $(1, 2n+2-j)$. In addition, the two extra classes also cross each other. So extra $(k-1)+(j-1)+1$ crossings are introduced.

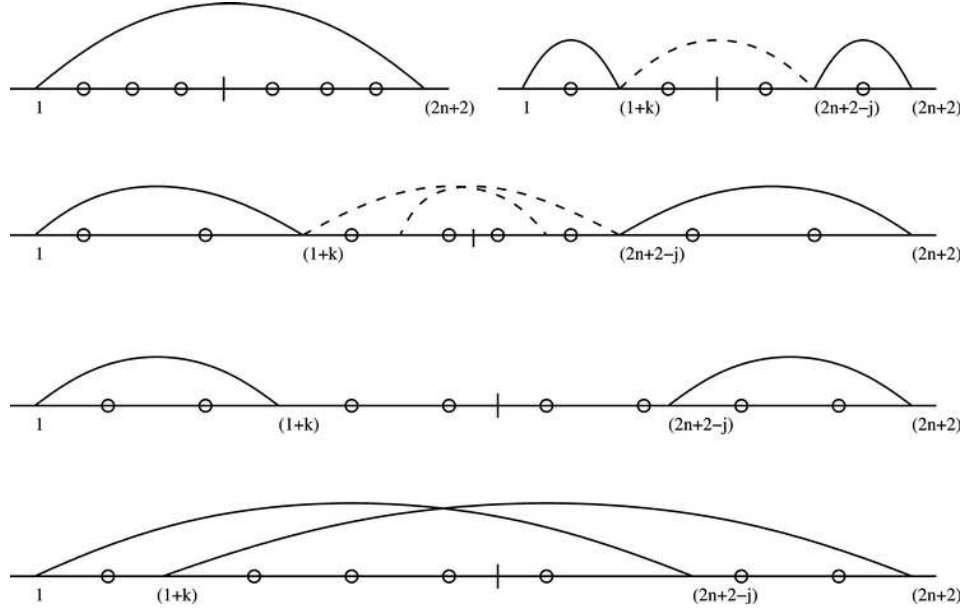

Fig. 3. *Five types of inhomogeneous partitions without singletons obtained from the product $\psi_1 \psi_n \psi_n \psi_1$.*

Taking expectations, using the fact that the process is centered, and the specific form (6) of the cumulants, we obtain

$$\mathrm{E}[\psi_1 \psi_n \psi_n \psi_1]$$
$$= \sum_{\pi \in \mathcal{P}_2(n,n)} q^{\mathrm{rc}(\pi)} t \cdot t^n$$

$$+ \sum_{\pi \in \mathcal{P}_2(n,n)} \left( \sum_{\substack{k,j=1,\ldots,n \\ k \sim j \mod \pi}} q^{k-1} q^{j-1} q^{\mathrm{rc}(\pi)} t^n + \sum_{\substack{k,j=1,\ldots,n \\ k \not\sim j \mod \pi}} q^{k-1} q^{j-1} q^{\mathrm{rc}(\pi)} t^n \right)$$

$$+ \sum_{k,j=1}^{n} \sum_{\pi \in \mathcal{P}_2(n-1,n-1)} q^{k-1} q^{j-1} q^{\mathrm{rc}(\pi)} t^2 t^{n-1}$$

$$+ \sum_{k,j=1}^{n} \sum_{\pi \in \mathcal{P}_2(n-1,n-1)} q \cdot q^{k-1} q^{j-1} q^{\mathrm{rc}(\pi)} t^2 t^{n-1}$$

$$= [n]_q! t^{n+1} + [n]_q^2 [n]_q! t^n + [n]_q^2 [n-1]_q! t^{n+1} + q[n]_q^2 [n-1]_q! t^{n+1}.$$

Similarly,

$$\mathrm{E}[\psi_{n+1} \psi_n X] = \mathrm{E}[X \psi_n \psi_{n+1}] = \mathrm{E}[\psi_{n+1} \psi_{n+1}] = [n+1]_q! t^{n+1},$$

$$\mathrm{E}[\psi_n \psi_n X] = \mathrm{E}[X \psi_n \psi_n] = [n]_q \mathrm{E}[\psi_n \psi_n] = [n]_q [n]_q! t^n$$

and

$$\mathrm{E}[\psi_{n-1} \psi_n X] = \mathrm{E}[X \psi_n \psi_{n-1}] = [n]_q t \mathrm{E}[\psi_{n-1} \psi_{n-1}] = [n]_q! t^n.$$

Substituting these relations into (4), we obtain

$$[n+1]_q! t^{n+1} + [n]_q^2 [n]_q! t^n + [n]_q [n]_q! t^{n+1}$$
$$+ [n]_q! t^{n+1} + [n]_q^2 [n]_q! t^n + (1+q)[n]_q [n]_q! t^{n+1}$$
$$- 2[n+1]_q! t^{n+1} - 2[n]_q^2 [n]_q! t^n - 2[n]_q [n]_q! t^{n+1}$$
$$= [n]_q! t^{n+1} + q[n]_q [n]_q! t^{n+1} - [n+1]_q! t^{n+1}$$
$$= (1 + q[n]_q - [n+1]_q)[n]_q! t^{n+1} = 0. \qquad \square$$

PROOF OF PROPOSITION 2.2. Let $Y \in \mathcal{A}_s$. Then

$$\mathrm{E}[Y \psi_n(t)] = \lim_{\delta(\mathcal{I}) \to 0} \mathrm{E}[Y \psi_n(t; \mathcal{I})].$$

Since the limit exists, we may restrict $\mathcal{I}$ to subdivisions containing $s$ as an endpoint of one of the intervals. The above expression is a sum of terms of the form

$$\mathrm{E}[Y X_{v(1)} X_{v(2)} \ldots X_{v(n)}].$$

If $I_{v(j)} \not\subset [0, s]$ for some $j$, the corresponding $X_{v(j)}$ is singleton independent from the rest of the terms in the product. Since the process is also centered, the resulting expectation is 0. As a result,

$$\mathrm{E}[Y \psi_n(t; \mathcal{I})] = \mathrm{E}[Y \psi_n(s; \mathcal{I})]$$

and so

$$\mathrm{E}[Y\psi_n(t)] = \mathrm{E}[Y\psi_n(s)].$$

Since this equality holds for an arbitrary $Y \in \mathcal{A}_s$, we conclude that the conditional expectation of $\psi_n(t)$ onto $\mathcal{A}_s$ is $\psi_n(s)$. $\square$

A transition operator for a Markov process is called Feller if it maps $C_0(\mathbb{R})$ into itself.

COROLLARY A.3. *Let $C_{q,n}$ be the scaled version of the continuous big $q$-Hermite polynomials and $\{X(t)\}$ be the centered $q$-Poisson process:*

(a) *$C_{q,n}(X(t), t)$ is a martingale with respect to the filtration induced by the process $\{X(t)\}$, for every $n$.*

(b) *Let*

$$H(x,t,z) = \prod_{k=0}^{\infty} \frac{1}{1 + ztq^k - (zq^k/(1+zq^k))(1-q)x}.$$

*Then $H(X(t), t, z)$ is a martingale.*

(c) *The process $\{X(t)\}$ is a Markov process with a Feller kernel.*

PROOF. Let $p_t$ be the orthogonal projection from $L^2(\mathbb{R}_+, dx)$ onto the subspace $L^2([0,t), dx)$. It can be extended to an operator on $\mathcal{F}_q(L^2(\mathbb{R}_+))$. The conditional expectation onto $\mathcal{A}_s$ is obtained by compression:

$$\mathrm{E}_t[A] = p_t A p_t.$$

The increments of a $q$-Lévy process are pyramidally, and so singleton, independent. Thus, the first part of the corollary follows from Propositions 2.2 and 5.5. It can also be obtained from the chaos decomposition property for the $q$-Poisson process,

$$C_{q,n}(X(t), t)\Omega = \mathbf{1}_{[0,t)}^{\otimes n}.$$

The second part follows from the first one since

$$H(x, t, z) = \sum_{n=0}^{\infty} \frac{1}{[n]_q!} C_{q,n}(x, t) z^n.$$

Note that the product defining $H$ converges for all $z$ since the sum

$$\sum_{k=0}^{\infty} \left( ztq^k - \frac{zq^k}{1+zq^k}(1-q)x \right)$$

converges.

LINEARIZATION COEFFICIENTS 23

The third part follows from the observations that the polynomials $\{C_{q,n}\}$ are, for every $t$, a basis for the polynomial ring, and polynomials are uniformly dense in the space of continuous functions on the (compact) spectrum of $X(t)$. Since the conditional expectation onto $\mathcal{A}_s$ is norm-continuous, this implies that for any continuous $f$, it maps $f(X_t)$ into the C*-algebra generated by $X(s)$. The existence of a Feller Markov kernel follows, see [5]. □

**Acknowledgments.** Thanks to Mourad Ismail and Dennis Stanton for useful discussions. Thanks also to the referees for a number of helpful comments and suggestions. Their diligent reading of the paper has resulted in a substantial improvement of the proofs.

DEPARTMENT OF MATHEMATICS
UNIVERSITY OF CALIFORNIA
RIVERSIDE, CALIFORNIA 92521-0135
USA
E-MAIL: manshel@math.ucr.edu